\documentclass[11pt]{article} 

\usepackage{algorithm, algorithmic, amsfonts, amsmath, amssymb, amsthm, graphicx, url}

\usepackage{fullpage}
\usepackage[small,bf]{caption}
\usepackage{color}
\usepackage[colorlinks,linkcolor=blue,citecolor=red]{hyperref}

\newcommand{\eq}[1]{\eqref{eq:#1}}
\newcommand{\eqd}[1]{\eqref{eq:#1}}
\newcommand{\fig}[1]{Figure~\ref{fig:#1}}

\numberwithin{equation}{section}

\newcommand{\R}{\mathbb{R}}  
\newcommand{\E}{\mathbb{E}}

\begin{document}

\title{A Tour of Reinforcement Learning:\\The View from Continuous Control}

\author{Benjamin Recht\\ 
Department of Electrical Engineering and Computer Sciences\\
University of California, Berkeley
}

\date{June 25, 2018. Last updated: November 9, 2018.}

\maketitle

\vspace{-0.25in}

\begin{abstract}

This manuscript surveys reinforcement learning from the perspective of optimization and control with a focus on continuous control applications. It surveys the general formulation, terminology, and typical experimental implementations of reinforcement learning and reviews competing solution paradigms.

In order to compare the relative merits of various techniques, this survey presents a case study of the Linear Quadratic Regulator (LQR) with unknown dynamics, perhaps the simplest and best-studied problem in optimal control. The manuscript describes how merging techniques from learning theory and control can provide non-asymptotic characterizations of LQR performance and shows that these characterizations tend to match experimental behavior. In turn, when revisiting more complex applications, many of the observed phenomena in LQR persist. In particular, theory and experiment demonstrate the role and importance of models and the cost of generality in reinforcement learning algorithms.

This survey concludes with a discussion of some of the challenges in designing learning systems that safely and reliably interact with complex and uncertain environments and how tools from reinforcement learning and control might be combined to approach these challenges.

\end{abstract}

\section{Introduction}

Reinforcement learning (RL) is the subfield of machine learning that studies how to use past data to enhance the future manipulation of a dynamical system. A control engineer might be puzzled by such a definition and interject that this is precisely the scope of control theory. That the RL and control communities remain practically disjoint has led to the co-development of vastly different approaches to the same problems. However, it should be impossible for a control engineer not to be impressed by the recent successes of the RL community such as solving Go~\cite{silver2016mastering}. 

Indeed, given this dramatic recent progress in reinforcement learning, a tremendous opportunity lies in deploying its data-driven systems in more demanding interactive tasks including self-driving vehicles, distributed sensor networks, and agile robotic systems.  For RL to expand into such technologies, however, the methods must be both safe and reliable---the failure of such systems has severe societal and economic consequences, including the loss of human life. How can we guarantee that our new data-driven automated systems are robust? These types of reliability concerns are at the core of  control engineering, and reinforcement learning practitioners might be able to make their methods robust by applying appropriate control tools for engineering systems to match prescribed safety guarantees.

This survey aims to provide a language for the control and reinforcement learning communities to begin communicating, highlighting what each can learn from the other. Controls is the theory of designing complex actions from well-specified models, while reinforcement learning often makes intricate, model-free predictions from data alone. Yet both RL and control aim to design systems that use richly structured perception, perform planning and control that adequately adapt to environmental changes, and exploit safeguards when surprised by a new scenario.  Understanding how to properly analyze, predict, and certify such systems requires insights from current machine learning practice and from the applied mathematics of optimization, statistics, and control theory.  With a focus on problems in continuous control, I will try to disentangle the similarities and differences of methods of the complementary perspectives and present a set of challenging problems whose solution will require significant input from both sets of practitioners.

I focus first on casting RL problems in an optimization framework, establishing the sorts of methodological tools brought to bear in contemporary RL. I then lay out the main solution techniques of RL including the dichotomy between the model-free and model-based methodologies. Next, I try to put RL and control techniques on the same footing through a case study of the linear quadratic regulator (LQR) with unknown dynamics. This baseline will illuminate the various trade-offs associated with techniques from RL and control. In particular, we will see that the so-called ``model-free'' methods popular in deep reinforcement learning are considerably less effective in both theory and practice than simple model-based schemes when applied to LQR. Perhaps surprisingly, I also show cases where these observations continue to hold on more challenging nonlinear applications. I then argue that model-free and model-based perspectives can be unified, combining their relative merits. This leads to a concluding discussion of some of the challenges at the interface of control and learning that must be solved before we can build robust, safe learning systems that interact with an uncertain physical environment, which will surely require tools from \emph{both} the machine learning and control communities.

\section{What is reinforcement learning?}

Reinforcement learning is the study of how to use past data to enhance the future manipulation of a dynamical system. How does this differ from ordinary machine learning? The main view of this survey is of reinforcement learning as optimal control when the dynamics are unknown. Our goal will be to find a sequence of inputs that drives a dynamical system to maximize some objective beginning with minimal knowledge of how the system responds to inputs.

In the classic optimal control problem, we begin with a dynamical system governed by the difference equation $x_{t+1} = f_t(x_t,u_t,e_t)$ where $x_t$ is the \emph{state} of the system, ${u_t}$ is the control action, and ${e_t}$ is a random disturbance. $f_t$ is the rule that maps the current state, control action, and disturbance at time $t$ to a new state. Assume that at every time, we receive some reward $R(x_t,u_t)$ for our current $x_t$ and $u_t$. The goal is to maximize this reward. In terms of mathematical optimization, we aim to solve the problem

\begin{equation}\label{eq:main-prob-tv}
\begin{array}{ll}
\mbox{maximize} & \E_{e_t}[ \sum_{t=0}^N R_t[x_t,u_t] ]\\
\mbox{subject to} &	x_{t+1} = f_t(x_t, u_t, e_t)\\
& \mbox{($x_0$ given).}
\end{array}
\end{equation}

That is, we aim to maximize the expected reward over $N$ time steps with respect to the control sequence $u_t$, subject to the dynamics specified by the state-transition rule $f_t$.  The expected value is over the disturbance, and assumes that $u_t$ is to be chosen having seen only the states $x_0$ through $x_t$ and previous inputs $u_0$ through $u_{t-1}$. $R_t$ is the reward gained at each time step and is determined by the state and control action. Note that $x_t$ is not really a decision variable in the optimization problem: it is determined entirely by the previous state, control action, and disturbance. I will refer to a \emph{trajectory}, $\tau_t$, as a sequence of states and control actions generated by a dynamical system.
\begin{equation}\label{eq:trajectory}
\tau_t = (u_1,\ldots,u_{t-1},x_0,\ldots,x_t) \,.
\end{equation}

Since the dynamics are stochastic, the optimal control problem typically allows a controller to observe the state before deciding upon the next action~\cite{BertsekasDPBook}. This allows a controller to continually mitigate uncertainty through \emph{feedback}. Hence, rather than optimizing over deterministic sequences of controls $u_t$, we instead optimize over $\emph{policies}$.  A \emph{control policy} (or simply ``a policy'') is a function, $\pi$, that takes a trajectory from a dynamical system and outputs a new control action.  Note that $\pi$  gets access only to previous states and control actions. 

To slightly lower the notational burden, I will hereon work with the \emph{time-invariant} version of Problem~\eqd{main-prob-tv}, assuming that the dynamical update rule is constant over time and that the rewards for state-action pairs are also constant:
\begin{equation}\label{eq:main-prob}
\begin{array}{ll}
\mbox{maximize} & \E_{e_t}[ \sum_{t=0}^N R(x_t,u_t) ]\\
\mbox{subject to} &	x_{t+1} = f(x_t, u_t, e_t),~u_t = \pi_t(\tau_t)\\
& \mbox{($x_0$ given).}
\end{array}
\end{equation}
The policies $\pi_t$ are the decision variables of the problem.

Let us now directly bring machine learning into the picture.  What happens when we don't know the state-transition rule $f$? There are a variety of commonly occurring scenarios when we might lack such knowledge. We may have unknown relationships between control forces and torques in a mechanical system. Or we could have a considerably more complicated system such as a massive data center with complex heat transfer interactions between the servers and the cooling systems. Can we still solve Problem~\eqd{main-prob} well without a precise model of the dynamics? Some lines of work even assume that we don't know the reward function $R$, but for the purpose of this survey, it makes no difference whether $R$ is known or unknown. The important point is that we can't solve this optimization problem using standard optimization methods unless we know the dynamics. We must learn something about the dynamical system and subsequently choose the best policy based on our knowledge.

The main paradigm in contemporary RL is to play the following game. We decide on a policy $\pi$ and horizon length L. Then we pass this policy either to a simulation engine or to a real physical system and are returned a trajectory $\tau_L$ and a sequence of rewards $\{R(x_t,u_t)\}$. We want to find a policy that maximizes the reward with the fewest total number of samples computed by the oracle, and we are allowed to do  whatever we'd like with the previously observed trajectories and reward information when computing a new policy. 
If we were to run $m$ queries with horizon length $L$, we would pay a total cost of $mL$.  However, we are free to vary our horizon length for each experiment. This is our \emph{oracle model} and is called \emph{episodic reinforcement learning} (See, for example Chapter 3 of Sutton and Barto~\cite{SuttonBartoBook}, Chapter 2 of Puterman~\cite{PutermanBook}, or Dann and Brunskill~\cite{Dann15}).  We want the expected reward to be high for our derived policy, but we also need the number of oracle queries to be small.

This oracle model is considerably more complicated than those typically considered in oracle models for optimization~\cite{NemirovskiYudinBook}. Each episode returns a complex feedback signal of states and rewards. What is the best way to tie this information together in order to improve performance? What is the best way to query and probe a system to achieve high quality control with as few interventions as possible? Here ``best'' is also not clearly defined. Do we decide an algorithm is best if it crosses some reward threshold in the fewest number of samples?  Or is it best if it achieves the highest reward given a fixed budget of samples? Or maybe there's a middle ground? This oracle provides a rich and complex model for interacting with a system and brings with it considerably more complexity than in standard stochastic optimization settings. What's the most efficient way to use all of the collected data in order to improve future performance?

\subsection{Connections to supervised learning}

The predominant paradigm of machine learning is \emph{supervised learning} or \emph{prediction}. In prediction, the goal is to predict the variable $y$ from a vector of \emph{features} $x$ such that on new data you are predicting $y$ from $x$ with high accuracy. This form of machine learning includes classification and regression as special cases. Most of the time when the term machine learning is used colloquially, it refers to this sort of prediction. From this perspective, niche topics like semi-supervised learning~\cite{Zhu-semi-sup} and matrix completion~\cite{hazan2012near} are prediction tasks as well.

By contrast,  there are \emph{two} special variables in reinforcement learning, $u$ and $r$. The goal now is to analyze the features $x$ and then subsequently choose a policy that emits $u$ so that $r$ is large.\footnote{To achieve notational consistency, I am throughout adopting the control-centric notation of denoting state-action pairs as $(x,u)$ rather than $(s,a)$ as is commonly used in reinforcement learning.} There are an endless number of problems where this formulation is applied~\cite{bertsekas1996neuro,kaelbling1996reinforcement,SuttonBartoBook} from online decision making in games~\cite{bowling2015heads,mnih2015human,silver2016mastering,tesauro1995td} to engagement maximization on internet platforms~\cite{bottou2013counterfactual,strehl2010learning}. A key distinguishing aspect of RL is the control action $u$. Unlike in prediction, the practitioner can vary $u$, which has implications both for learning (e.g., designing experiments to learn about a given system) and for control (e.g., choosing inputs to maximize reward).

Reinforcement learning is clearly more challenging than supervised learning, but, at the same time, it can be considerably more valuable. Reinforcement learning provides a useful framework to conceptualize interaction in machine learning, and promises to help mitigate changing distributions, gaming, adversarial behavior, and unexpected amplification. There is a precarious trade-off that must be carefully considered:  reinforcement learning demands interventions with the promise that these actions will directly lead to valuable returns, but the resulting complicated feedback loops are hard to study in theory, and failures can have catastrophic consequences.

\section{Strategies for solving reinforcement learning problems}

Let us now turn to a taxonomy of the varied algorithmic frameworks for reinforcement learning, focused on solving Problem~\eqd{main-prob} when the state-transition function is unknown. \emph{Model-Based} Reinforcement learning fits a model to previously observed data and then uses this model in some fashion to approximate the solution to Problem~\eqd{main-prob}. \emph{Model-Free} Reinforcement learning eschews the need for a system's model, directly seeking a map from observations to actions.  

The role of models in reinforcement learning remains hotly debated. Model-free methods, as discussed below, aim to solve optimal control problems only by probing the system and improving strategies based on past rewards and states. Many researchers argue for algorithms that can innately learn to control without access to the complex details required to simulate a dynamical system. They argue that it is often easier to find a policy for a task than it is to fit a general purpose model of the system dynamics (see for example, the discussion in Chapter 3 of Volume 2 of Bertsekas~\cite{BertsekasDPBook2}). Model-free methods are primarily divided into two approaches: \emph{Policy Search} and \emph{Approximate Dynamic Programming}. Policy Search directly searches for policies by using data from previous episodes in order to improve the reward. Approximate Dynamic Programming uses Bellman's principle of optimality to approximate Problem~\eqd{main-prob} using previously observed data. 

Throughout, my aim will be to highlight the main conceptual ideas of different approaches and to avoid embroiling myself in a thorough discussion of the myriad of technical details required to make all of the statements crisply precise. What is important is that all of the approaches surveyed reduce to some sort of function fitting from noisy observations of the dynamical system, though performance can be drastically different depending on how you fit this function. In model-based reinforcement learning, we fit a model of the state transitions to best match observed trajectories. In approximate dynamic programming, we estimate a function that best characterizes the ``cost to go'' for experimentally observed states. And in direct policy search, we attempt to find a policy that directly maximizes the optimal control problem using only input-output data. The main question is which of these approaches makes the best use of samples and how quickly do the derived policies converge to optimality. 

\subsection{Model-based reinforcement learning} 

One of the simplest and perhaps most obvious strategies to solve the core RL Problem~\eqd{main-prob} is to estimate a predictive model for the dynamical process and then to use it in a dynamic programming solution to the prescribed control problem. The estimated model is called the \emph{nominal model}, and I will refer to control design assuming the estimated model is true as \emph{nominal control}. Nominal control, commonly verbosely referred to as ``control under the principle of certainty equivalence,'' serves as a useful baseline algorithm. 

Estimation of dynamical systems is called  \emph{system identification} in the control community~\cite{LjungBook}. System Identification differs from conventional estimation because one needs to carefully choose the right inputs to excite  various degrees of freedom and because dynamical outputs are correlated over time with the parameters we hope to estimate, the inputs we feed to the system, and the stochastic disturbances. Once data are collected, however, conventional machine learning tools can be used to find the system that best agrees with the data and can be applied to analyze the number of samples required to yield accurate models~\cite{campi2002finite,vidyasagar2008learning}.

Suppose we want to build a predictor of $x_{t+1}$ from the trajectory history. A simple, classic strategy is simply to inject a random probing sequence $u_t$ for control and then measure how the state responds.  Up to stochastic noise, we should have that
\begin{equation}\label{eq:sys-id}
	x_{t+1} \approx \varphi(x_t,u_t) \,,
\end{equation}
where $\varphi$ is some model aiming to approximate the true dynamics. $\varphi$ might arise from a first-principles physical model or might be a non-parametric approximation by a neural network. The state-transition function can then be fit using supervised learning. For instance, a model can be fit by solving the least squares problem
\begin{equation*}
\begin{array}{ll}
\mbox{minimize}_{\varphi} & \sum_{t=0}^{N-1} ||x_{t+1} - \varphi(x_t,u_t)||^2\,.
\end{array}
\end{equation*}

Let $\hat{\varphi}$ denote the function fit to the collected data to model the dynamics. Let $\omega_t$ denote a random variable that we will use as a model for the noise process. With such a point estimate for the model, we might solve the optimal control problem
\begin{equation*}\label{eq:main-prob-approx-dyn}
\begin{array}{ll}
\mbox{maximize} & \E_{\omega_t}[ \sum_{t=0}^N R(x_t,u_t)  ]\\
\mbox{subject to} &	x_{t+1} = \hat{\varphi}(x_t, u_t)+\omega_t,~u_t = \pi_t(\tau_t)\,.
\end{array}
\end{equation*}
In this case, we are solving the wrong problem to get our control policies $\pi_t$. Not only is the model incorrect, but this formulation requires some plausible model of the noise process. But if $\hat{\varphi}$ and $f$ are close, this approach might work well in practice. 

\subsection{Approximate Dynamic Programming}\label{sec:adp}

Approximate dynamic programming approaches the RL problem by directly approximating the optimal control cost and then solving this approximation with techniques from dynamic programming. The dynamic programming solution to Problem~\eqd{main-prob} is based on the \emph{principle of optimality}: if you've found an optimal control policy for a time horizon of length $N$, $\pi_1,\ldots, \pi_N$, and you want to know the optimal strategy starting at state $x$ at time $t$, then you just have to take the optimal policy starting at time $t$, $\pi_t,\ldots,\pi_N$. Dynamic programming then lets us recursively find a control policy by starting at the final time and recursively solving for policies at earlier times.

Define the \emph{Q-function} for~\eq{main-prob} to be the mapping 
\begin{equation}\label{eq:q-fun}
	\mathcal{Q}(x,u) = \max\left\{ \E_{e_t}\left[ \sum_{t=0}^N R(x_t,u_t)\right] \,:\,x_{t+1} = f(x_t, u_t, e_t),\,(x_0,u_0)=(x,u)\right\}
\end{equation}
The Q-function determines the value of the optimal control problem that is attained when the first action is set to be $u$ and the initial condition is $x$. Note that it then trivially follows that the optimal value of Problem~\eqd{main-prob} is $\max_u \mathcal{Q}(x_0,u)$, and the optimal policy is $\pi(x_0) = \arg \max_u \mathcal{Q}(x_0,u)$. If we had access to the Q-function, we'd have everything we'd need to know to take the first step in the optimal control problem. We can use dynamic programming to compute this Q-function and the Q-function associated with every subsequent action. That is, we define the terminal Q-function to be
\begin{equation*}
	\mathcal{Q}_{N}(x,u) = R(x,u)\,,
\end{equation*}
and then define recursively
\begin{equation}\label{eq:bellman-finite}
	\mathcal{Q}_{k}(x,u) =R(x,u) + \E_{e}\left[ \max_{u'} \mathcal{Q}_{k+1}(f(x,u,e),u')\right]\,.
\end{equation}
This is the dynamic programing algorithm in a nutshell: we can recursively define the Q-functions by passing backward in time, and then compute the optimal controls from any starting $x_0$ by applying the policy that maximizes the right hand side of~\eq{bellman-finite} at each time step.~\eq{bellman-finite} is known as Bellman's equation. Note that for all time, the optimal policy is $u_k = \arg\max_u \mathcal{Q}_k (x_k,u)$ and \emph{depends only on the current state.}

Approximate Dynamic Programming methods typically try to compute these action-value functions from data. They do so by assuming that the Q-function is stationary. (i.e., $\mathcal{Q}_k(x,u) = \mathcal{Q}(x,u)$ for all $k$ and some function $\mathcal{Q}$). Such stationarity indeed arises assuming the time horizon is infinite. Consider the limit:
\begin{equation}\label{eq:main-prob-inf}
\begin{array}{ll}
\mbox{maximize} & \lim_{N\rightarrow \infty}  \E_{e_t}[ \frac{1}{N} \sum_{t=0}^N R(x_t,u_t) ]\\
\mbox{subject to} &	x_{t+1} = f(x_t, u_t, e_t),~u_t=\pi_t(\tau_t)\\
& \mbox{($x_0$ given).}
\end{array}
\end{equation}
And we define the Q-function $\mathcal{Q}(x_0,u_0)$ to be the average reward accrued running from state $x_0$ with initial action $u_0$.
Unfortunately, Problem~\eqd{main-prob-inf} is not directly amenable to dynamic programming without introducing further technicalities. For
mathematical convenience and also to connect to common practice in RL, it's useful to instead consider the \emph{discounted} reward problem
\begin{equation}\label{eq:main-prob-disc}	
	\begin{array}{ll}	
		\mbox{maximize} &  (1-\gamma) \E_{e_t}[ \sum_{t=0}^\infty \gamma^t R(x_t,u_t) ]\\
		\mbox{subject to} &	x_{t+1} = f(x_t, u_t, e_t),~u_t=\pi_t(\tau_t)\\
		& \mbox{($x_0$ given).}
	\end{array}
\end{equation}
where $\gamma$ is a scalar in $(0,1)$ called the \emph{discount factor}. For $\gamma$ close to $1$, the discounted reward is approximately equal to the average reward~\cite{BertsekasDPBook2}. The discounted cost has particularly clean optimality conditions that make it amenable to estimation.
If we define $\mathcal{Q}_\gamma(x,u)$ to be the Q-function obtained from solving Problem~\eqd{main-prob-disc} with initial condition $x$, then we have a discounted version of dynamic programming, now with the same Q-functions on the left and right hand sides:
\begin{equation*}\label{eq:bellman}
	\mathcal{Q}_\gamma(x,u) = R(x,u) + \gamma \E_{e} \left[ \max_{u'} \mathcal{Q}_\gamma(f(x,u,e),u')\right]\,.
\end{equation*}
The optimal policy is now for \emph{all times} to let 
\begin{equation}\label{eq:opt-q-policy}
	u_t= \arg\max_u \mathcal{Q}_\gamma(x_t,u)\,.
\end{equation} 
This is a remarkably simple formula which is part of what makes Q-learning methods so attractive.

We can try to solve for the Q-function using stochastic approximation. If we draw a sample trajectory using the policy given by~\eq{opt-q-policy}, then we should have (approximately and in expectation)
\[
	\mathcal{Q}_\gamma(x_k,u_k) \approx R(x_k,u_k) + \gamma \max_{u'} \mathcal{Q}_\gamma(x_{k+1},u')\,.
\]
Thus, beginning with some initial guess $\mathcal{Q}_\gamma^{(\mathrm{old})}$ for the Q-function, we can update
\begin{equation}\label{eq:q-learning}
	Q_\gamma^{(\mathrm{new})}(x_k,u_k) = (1-\eta) Q_\gamma^{(\mathrm{old})}(x_k,u_k) + \eta \left(R(x_k,u_k) + \gamma \max_{u'}  Q_\gamma^{(\mathrm{old})}(x_{k+1},u')\right)
\end{equation}
where $\eta$ is a \emph{step-size} or \emph{learning rate}.~\eq{q-learning} forms the basis of \emph{Q-learning} algorithms~\cite{tsitsiklis1994asynchronous,watkins1992q}.

Surveying ADP using only Q-functions is somewhat unorthodox. Most introductions to ADP instead  focus on \emph{value} functions where
\[
	V(x) = \max_u \mathcal{Q}(x,u)\,.
\]
Methods for estimating value functions are also widely used in reinforcement learning and developed through the perspective of estimation and stochastic approximation. In particular, Temporal Difference algorithms are derived from the value-function-centric perspective~\cite{sutton1988learning,dayan1992convergence,bradtke1996linear,bertsekas1996temporal,yu2009convergence}. 

Note that in all cases here, though we have switched away from models, there's no free lunch. We are still estimating functions here, and we need to assume that the functions have some reasonable structure or we can't learn them. Choosing a parameterization of the Q-function \emph{is a modeling assumption.} The term ``model-free'' almost always means ``no model of the state transition function'' when casually claimed in reinforcement learning research. However, this does not mean that modeling is not heavily built into the assumptions of model-free RL algorithms. Moreover, for continuous control problems these methods appear to make an inefficient use of samples. Suppose the internal state of the system is of dimension $d$. When modeling the state-transition function,~\eq{sys-id} provides $d$ equations per time step. By contrast, we are only using $1$ equation per time step in ADP. Such inefficiency is certainly seen in practice below. Also troubling is the fact that we had to introduce the discount factor in order to get a simple Bellman equation. One can avoid discount factors, but this requires considerably more sophisticated analysis. Large discount factors do in practice lead to brittle methods, and the discount becomes a hyperparameter that must be tuned to stabilize performance. We will determine below when and how these issues arise in practice in control.

\subsection{Direct Policy Search}\label{sec:policy-search}

The most ambitious form of control without models attempts to directly learn a policy function from episodic experiences without ever building a model or appealing to the Bellman equation. From the oracle perspective, these policy driven methods turn the problem of RL into derivative-free optimization.

\subsubsection{A generic algorithm for sampling to optimize}

In turn, let's first begin with a review of a general paradigm for leveraging random sampling to solve optimization problems. Consider the general unconstrained optimization problem
\begin{equation}\label{eq:generic-reward-opt}
\begin{array}{ll}
	\mbox{maximize}_{z\in\R^d} & R(z) \,.
	\end{array}
\end{equation}
Any optimization problem like this is equivalent to an optimization over probability distributions on $z$:
$$
\begin{array}{ll}
	\mbox{maximize}_{p(z)} & \mathbb{E}_p[R(z)] \,.
\end{array}
$$
If $z_\star$ is the optimal solution, then we'll get the same value if we put a $\delta$-function around $z_\star$.  Moreover, if $p$ is a probability distribution, it is clear that the \emph{expected value of the reward function} can never be larger than the maximal reward achievable by a fixed $z$. So we can either optimize over $z$ or we can optimize over \emph{distributions} over $z$.

Since optimizing over the space of all probability densities is intractable, we must restrict the class of densities over which we optimize. For example, we can consider a family parameterized by a parameter vector $\vartheta$: $p(u;\vartheta)$ and attempt to optimize
\begin{equation}\label{eq:stoch-search}
\begin{array}{ll}
	\mbox{maximize}_{\vartheta} & \mathbb{E}_{p(z;\vartheta)}[R(z)] \,.
\end{array}
\end{equation}
If this family of distributions contains all of the Delta functions, then the optimal value will coincide with the non-random optimization problem. But if the family does not contain the Delta functions, the resulting optimization problem only provides a lower bound on the optimal value no matter how good of a probability distribution we find. 

That said, this reparameterization provides a powerful and general algorithmic framework for optimization. In particular, we can compute the derivative of $J(\vartheta):= \mathbb{E}_{p(z;\vartheta)}[R(z)] $ using the following calculation (called ``the log-likelihood trick''):
\begin{align*}
	\nabla_{\vartheta} J(\vartheta) &= \int R(z) \nabla_{\vartheta} p(z;\vartheta) dz\\
	&= \int R(z) \left(\frac{\nabla_{\vartheta} p(z;\vartheta)}{p(z;\vartheta)}\right) p(z;\vartheta) dz\\
	&= \int \left( R(z) \nabla_{\vartheta} \log p(z;\vartheta) \right) p(z;\vartheta)dz	\\
  &= \mathbb{E}_{p(z;\vartheta)}\left[ R(z) \nabla_{\vartheta} \log p(z;\vartheta) \right]\,.
\end{align*}
This derivation reveals that the gradient of $J$ with respect to $\vartheta$ is the expected value of the function
\begin{equation}\label{eq:policy-grad}
	G(z,\vartheta) = R(z) \nabla_{\vartheta} \log p(z;\vartheta)
\end{equation}
Hence, if we sample $z$ from the distribution defined by $p(z;\vartheta)$, we can compute $G(z,\vartheta)$ and will have an unbiased estimate of the gradient of $J$. We can follow this direction and will be running stochastic gradient descent on $J$, defining Algorithm~\ref{alg:reinforce}.

\begin{center}
\begin{algorithm}[ht]
\begin{algorithmic}[1]
\STATE {\bf Hyperparameters:} step-sizes $\alpha_j>0$.
\STATE {\bf Initialize:} $\vartheta_0$ and $k = 0$.
\WHILE{ending condition not satisfied}
\STATE Sample $z_k\sim p(z;\vartheta_k)$.
\STATE Set $\vartheta_{k+1} = \vartheta_k + \alpha_k R(z_k) \nabla_\vartheta \log p(z_k; \vartheta_k)$.
\STATE $k\leftarrow k + 1$
\ENDWHILE
\end{algorithmic}
\caption{REINFORCE}
\label{alg:reinforce}
 \end{algorithm}
\end{center} 

Algorithm~\ref{alg:reinforce} is typically called REINFORCE~\cite{williams1992simple} and its main appeal is that it is trivial to implement. If you can efficiently sample from $p(z;\vartheta)$, you can run this algorithm on essentially any problem. But such generality must and does come with a significant cost. The algorithm operates on stochastic gradients of the sampling distribution, but the function we cared about optimizing---$R$---is only accessed through function evaluations. Direct search methods that use the log-likelihood trick are necessarily derivative free optimization methods, and, in turn, are necessarily less effective than methods that compute actual gradients, especially when the function evaluations are noisy~\cite{Jamieson12}. Another significant concern is that the choice of distribution can lead to very high variance in the stochastic gradients. Such high variance in turn implies that many samples need to be drawn to find a stationary point.

That said, the ease of implementation should not be readily discounted. Direct search methods are trivial to implement, and oftentimes reasonable results can be achieved with considerably less effort than custom solvers tailored to the structure of the optimization problem. There are two primary ways that this sort of stochastic search arises in reinforcement learning: Policy Gradient and Pure Random Search.

\subsubsection{Policy Gradient}

As seen from Bellman's equation, the optimal policy for Problem~\eqd{main-prob} is always deterministic. Nonetheless, the main idea behind policy gradient is to use \emph{probabilistic policies}. Probabilistic policies are optimal for other optimization-based control problems such as control of partially observed Markov decision processes~\cite{aastrom1965optimal,kaelbling1998planning} or in zero-sum games. Hence, exploring their value for the RL problems studied in this survey does not appear too outlandish at first glance.

We fix our attention on \emph{parametric, randomized policies} such that $u_t$ is sampled from a distribution $p(u \vert \tau_t;\vartheta)$ that is a function only of the currently observed trajectory and a parameter vector $\vartheta$. A probabilistic policy induces a probability distribution over trajectories:
\begin{equation}\label{eq:factor-graph}
	p(\tau;\vartheta) = \prod_{t=0}^{L-1} p(x_{t+1} \vert x_{t},u_{t}) p(u_t \vert \tau_t ;\vartheta)\,.
\end{equation}
Moreover, we can overload notation and define the reward of a trajectory to be
\begin{equation*}
	R(\tau) = \sum_{t=0}^N R_t(x_t,u_t)
\end{equation*}
Then our optimization problem for reinforcement learning tidily takes the form of Problem~\eqd{stoch-search}. Policy gradient thus proceeds by sampling a trajectory using the probabilistic policy with parameters $\vartheta_k$, and then updating using REINFORCE.

Using the log-likelihood trick and~\eq{factor-graph}, it is straightforward to verify that the gradient of $J$ with respect to $\vartheta$ \emph{is not an explicit function of the underlying dynamics}.  However, at this point this should not be surprising. By shifting to distributions over policies, we push the burden of optimization onto the sampling procedure. 
\subsubsection{Pure Random Search}

An older and more widely applied method to solve Problem~\eqd{generic-reward-opt} is to directly perturb the current decision variable $z$ by random noise and then update the model based on the received reward at this perturbed value. That is, we apply Algorithm~\ref{alg:reinforce} with sampling distribution $p(z;\vartheta) = p_0(z-\vartheta)$ for some distribution $p_0$. The simplest examples for $p_0$ would be the uniform distribution on a sphere or a normal distribution. Perhaps less surprisingly here, REINFORCE can again be run without any knowledge of the underlying dynamics. Note that in this case, the REINFORCE algorithm has a simple interpretation in terms of gradient approximation. Indeed, REINFORCE is equivalent to approximate gradient ascent of $R$
$$
	\vartheta_{t+1} = \vartheta_{t} + \alpha g_\sigma(\vartheta_k)
$$
 with the gradient approximation
$$
	g_\sigma(\vartheta) = \frac{R(\vartheta + \sigma \epsilon) - R(\vartheta - \sigma \epsilon) }{2\sigma} \epsilon\,.
$$
This update says to compute a finite difference approximation to the gradient along the direction $\epsilon$ and move along the gradient. One can reduce the variance of such a finite-difference estimate by sampling along multiple random directions and averaging:
$$
	g^{(m)}_\sigma(\vartheta) = \frac{1}{m} \sum_{i=1}^m\frac{R(\vartheta + \sigma \epsilon_i) - R(\vartheta - \sigma \epsilon_i) }{2\sigma} \epsilon_i\,.
$$
This is akin to approximating the gradient in the random subspace spanned by the $\epsilon_i$

This particular algorithm and its generalizations go by many different names. Probably the earliest proposal for this method was made by Rastrigin~\cite{Rastrigin63}. In an unexpected historical surprise, Rastrigin initially developed this method to solve reinforcement learning problems! His main motivating example was an inverted pendulum. A rigorous analysis using contemporary techniques was provided by Nesterov and Spokoiny~\cite{nesterov2017random}. Random search was also discovered by the evolutionary algorithms community, where it is called a $(\mu,\lambda)$-Evolution Strategy~\cite{Beyer02,SchwefelThesis}. Random search has also been studied in the context of stochastic approximation~\cite{spall1992multivariate} and bandits~\cite{flaxman2005online,agarwal2010optimal}. Algorithms that are invented independently by four different communities probably have something good going for them.

The random search method is considerably simpler than the policy gradient algorithm but it uses much less structure from the problem as well.  Since RL problems tend to be nonconvex, it is not clear which of these approaches is best unless we focus on specific instances. In light of this, in the next section we turn to a set of instances where we may be able to glean more insights about the relative merits of all of the approaches to RL covered in this section.

\subsection{Deep reinforcement learning}

Note that in this section I have spent no time discussing \emph{deep} reinforcement learning. That is because there is nothing conceptually different other than the use of neural networks for function approximation. That is, if one wants to take any of the described methods and make them deep, they simply need to add a neural net. In model-based RL, $\varphi$ is parameterized as a neural net, in ADP, the Q-functions or Value Functions are assumed to be well-approximated by neural nets, and in policy search, the policies are set to be neural nets. The algorithmic concepts themselves don't change. However, convergence analysis certainly will change, and algorithms like Q-learning might not even converge. The classic text Neuro-dynamic Programming by Bertsekas and Tsitisklis discusses the adaptations needed to admit function approximation~\cite{bertsekas1996neuro}. By eliminating the complicating variable of function approximation, we can get better insights into the relative merits of these methods, especially when focusing on a simple set of instances of optimal control, namely, the Linear Quadratic Regulator.

\section{Simplifying theme: The Linear Quadratic Regulator}

With this varied list of approaches to reinforcement learning, it is difficult from afar to judge which method fares better on which problems. It is likely best to start simple and small and find the simplest non-trivial problem that can assist in distinguishing the various approaches to control.  Though simple models are not the end of the story in analysis, it tends to be the case that if a complicated method fails to perform on a simple problem, then this indicates a flaw in the method.

I'd argue that in controls, the simplest non-trivial class of instances of optimal control is those with convex quadratic rewards and linear dynamics. That is, the problem of the Linear Quadratic Regulator (LQR):
\begin{equation}\label{eq:lqr-prob}
\begin{array}{ll}
\mbox{minimize} \, & \E_{e_t} \left[\frac{1}{2}\sum_{t=0}^N x_t^TQ x_t + u_t^T R u_t  + \frac{1}{2} x_{N+1}^T S x_{N+1}\right], \\
\mbox{subject to} & x_{t+1} = A x_t+ B u_t +e_t,~u_t=\pi_t(\tau_t) \\
& \mbox{($x_0$ given).}
\end{array}
\end{equation}
Here, $Q$, $R$, and $S$ are positive semidefinite matrices. Do note that we have switched to minimization from maximization, as is conventional in optimal control. The state transitions are governed by a linear update rule with $A$ and $B$ appropriately sized matrices.

A few words are in order to defend this baseline as instructive for general problems in continuous control and RL. Though linear dynamics are somewhat restrictive, many systems are linear over the range we'd like them to operate. Indeed, enormous engineering effort goes into designing systems so that their responses are as close to linear as possible. From an optimization perspective, linear dynamics are the only class where we are guaranteed that our constraint set is convex, which is another appealing feature for analysis.

What about cost functions? Whereas dynamics are typically handed to the engineer, cost functions are completely at their discretion. Designing and refining cost functions are part of optimal control design, and different characteristics can be extracted by iteratively refining cost functions to meet specifications. This is no different in machine learning where, for example, combinatorial losses in classification are replaced with smooth losses like logistic or squared loss. Designing cost functions is a major challenge and tends to be an art form in engineering. But since we're designing our cost functions, we should focus our attention on costs that are easier to solve. Quadratic cost is particularly attractive not only because it is convex, but also for how it interacts with noise. The cost of the stochastic problem is equal to that of the noiseless problem plus a constant that is independent of the choice of $u_t$. The noise will degrade the achievable cost, but it will not affect how control actions are chosen. 

Note that when the parameters of the dynamical system are known, the standard LQR problem admits an elegant dynamic programming solution~\cite{Zhou95}. The control action is a \emph{linear function} of the state
\begin{equation*}
	u_t = -K_t x_t
\end{equation*}
for some matrix $K_t$ that can be computed via a simple linear algebraic recursion with only knowledge of $(A,B,Q,R)$.

In the limit as the time horizon tends to infinity, the optimal control \emph{policy} is \emph{static}, \emph{linear} \emph{state feedback}:
\begin{equation*}
	u_t = -K x_t
\end{equation*}
where $K$ is a fixed matrix defined by
\begin{equation*}
	K=(R + B^T M B)^{-1} B^T M A
\end{equation*}
and $M$ is a solution to the \emph{Discrete Algebraic Riccati Equation}
\begin{equation}\label{eq:dare}
M = Q + A^T M A - (A^T M B)(R + B^T M B)^{-1} (B^T M A)\,.
\end{equation}
That is, for LQR on an infinite time horizon, $\pi_t (x_t) = - K x_t$. Here, $M$ is the unique solution of the Riccati equation where
all of the eigenvalues of $A-BK$ have magnitude less than $1$. Finding this specific solution
is relatively easy using standard linear algebraic techniques~\cite{Zhou95}.

There are a variety of ways to derive these formulae. In particular, one can use dynamic programming as in Section~\ref{sec:adp}. In this case, one can check that the Q-function on a finite time horizon satisfies a recursion
\begin{equation*}
 \mathcal{Q}_k(x,u) = x^T Q x + u^T R u + (Ax+Bu)^T M_{k+1} (Ax+Bu) + c_k\,.
\end{equation*}
for some positive definite matrix $M_{k+1}$. The limit of these matrices are the solution of~\eq{dare}

Though LQR cannot capture every interesting optimal control problem, it has many of the salient features of the generic optimal control problem. Dynamic programming recursion lets us compute the control actions efficiently and, for long time horizons, a static policy is nearly optimal.

Now the main question to consider in the context of RL: What happens when we don't know $A$ and $B$? What's the right way to interact with the dynamical system in order to quickly and efficiently get it under control? Let us now dive into the different styles of reinforcement learning and connect them to ideas in controls, using LQR as a guiding baseline.

\subsection{The sample complexity of model-based RL for LQR}\label{sec:nominal}

For LQR, maximum likelihood estimation of a nominal model is a least squares problem:
$$
\begin{array}{ll}
\mbox{minimize}_{A,B} & \sum_{t=0}^{N-1} ||x_{t+1} - A x_t -  B u_t||^2\,.
\end{array}
$$
How well do these model estimates work for the LQR problem? Suppose we treat the estimates as true and use them to compute a state feedback control from a Riccati equation. While we might expect this to work well in practice, how can we verify the performance?  As a simple case, suppose that the true dynamics are slightly unstable so that $A$ has at least one eigenvalue of magnitude larger than $1$.  It is fully possible for the least squares estimates of such a mode is less than one, and, consequently, the optimal control strategy using the estimate will fail to account for the poorly estimated unstable eigenvalue. How can we include the knowledge that our model is just an estimate and not accurate with a small sample count? One possible solution is to use tools from robust control to mitigate this uncertainty.

\subsubsection{Coarse-ID Control: a new paradigm for learning to control.}

My collaborators and I have been considering an approach to merge robust control and high-dimensional statistics dubbed ``Coarse-ID Control.'' The general framework consists of the following three steps:  
\begin{enumerate}
\item Use supervised learning to learn a coarse model of the dynamical system to be controlled.  I will refer to the system estimate as the \emph{nominal system}.
\item Using either prior knowledge or statistical tools like the bootstrap, build probabilistic guarantees about the distance between the nominal system and the true, unknown dynamics.
\item Solve a \emph{robust optimization} problem that optimizes control of the nominal system while penalizing signals with respect to the estimated uncertainty, ensuring stable, robust execution.
\end{enumerate}

As long as the true system behavior lies in the estimated uncertainty set, we'll be guaranteed to find a performant controller.  The key here is that we are using machine learning to identify not only the plant to be controlled, \emph{but the uncertainty as well}. Indeed, the main advances in the past two decades of estimation theory consist of providing reasonable estimates of such uncertainty sets with guaranteed bounds on their errors as a function of the number of observed samples. Taking these new tools and merging them with old and new ideas from robust control allow us to bound the end-to-end performance of a controller in terms of the number of observations.

The coarse-ID procedure is well illustrated through the case study of LQR~\cite{Dean17}. We can guarantee the accuracy of the least squares estimates for $A$ and $B$ using novel probabilistic analysis~\cite{Simchowitz18a}. With the estimate of model error in hand, one can pose a robust variant of the standard LQR optimal control problem that computes a robustly stabilizing controller seeking to minimize the worst-case performance of the system given the (high-probability) norm bounds on our modeling errors.  

To design a good control policy, we here turn to state-of-the-art tools from robust control. We leverage the recently developed System Level Synthesis (SLS) framework~\cite{Matni17,Wang16} to solve this robust optimization problem. SLS lifts the system description into a higher dimensional space that enables efficient search for controllers.  The proposed approach provides non-asymptotic bounds that guarantee finite performance on the infinite time horizon, and quantitatively bound the gap between the computed solution and the true optimal controller.

Suppose in LQR that we have a state dimension $d$ and control dimension $p$. Denote the minimum cost achievable by the optimal controller as $J_\star$. Our analysis guarantees that after a observing a trajectory of length $T$, we can design a controller that will have infinite-time-horizon cost $\hat{J}$ with 
\begin{equation*}
\frac{\hat{J}-J_\star}{J_\star} = \tilde{O}\left(\sqrt{\tfrac{d+p}{T}}\right)\,.
\end{equation*}
Here, the notation $\tilde{O}(\cdot)$ suppresses logarithmic factors and instance-dependent constants.  In particular, we can guarantee that we stabilize the system after seeing only a finite amount of data.

Where Coarse-ID control differs from nominal control is that it explicitly accounts for the uncertainty in the least squares estimate. By appending this uncertainty to the original LQR optimization problem, we can circumvent the need to study perturbations of Riccati equations. Moreover, since the approach is optimization based, it can be readily applied to other optimal control problems beyond the LQR baseline.

\subsection{The sample complexity of model-free RL for LQR}

Since we know that the Q-function for LQR is quadratic, we can try to estimate it by dynamic programming. Such a method was probably first proposed by Bradtke, Barto, and Ydstie~\cite{bradtke1994adaptive}. More recently, Tu showed that the \emph{Least-squares Temporal Differencing} algorithm, also due to Bradtke and Barto~\cite{bradtke1996linear}, could estimate the value function of LQR to low error with $\tilde{O}\left(\sqrt{\tfrac{d^2}{T}}\right)$ samples~\cite{Tu17b}. This estimator can then be combined with a method to improve the estimated policy over time.

Note that the bound on the efficiency of the estimator here is worse than the error obtained for estimating the model of the dynamical system. While comparing worst-case upper bounds is certainly not valid, it is suggestive that, as mentioned above, temporal differencing methods use only one defining equation per time step whereas model estimation uses $d$ equations per time step. So while the conventional wisdom suggests that estimating Q-functions for specific tasks should be simpler than estimating models, the current methods appear to be less efficient with aggregated data than system identification methods.

With regard to direct search methods, we can already see variance issues enter the picture even for small LQR instances. Consider the most trivial example of LQR:
\begin{equation*}
	R(u) = ||u||^2
\end{equation*}
Let $p(u;\vartheta)$ be a multivariate Gaussian with mean $\vartheta$ and variance $\sigma^2 I$.  Then
\begin{equation*}
	\mathbb{E}_{p(u;\vartheta)} [R(u)]= \|\vartheta\|^2 + \sigma^2 d
\end{equation*}
Obviously, the best thing to do would be to set $\vartheta=0$. Note that the expected reward is off by $\sigma^2 d$ at this point, but at least this would be finding a good guess for $u$.  Also, as a function of $\vartheta$, the cost is \emph{strongly convex}, and the most important thing to know is the expected norm of the gradient as this will control the number of iterations. Now, after sampling $u$ from a Gaussian with mean $\vartheta_0$ and variance $\sigma^2 I$ and using formula~\eqd{policy-grad}, the first gradient will be  
\begin{equation*}
	g=-\frac{||\omega-\vartheta_0||^2 \omega}{\sigma^2}\,,
\end{equation*}
where $\omega$ is a normally distributed random vector with mean zero and covariance $\sigma^2 I$.
The expected norm of this stochastic gradient is on the order of
\begin{equation*}
O\left(\sigma d^{1.5} + \sigma^{-1} d^{0.5} \|\vartheta_0\|\right)\,,
\end{equation*}
which indicates a significant scaling with dimension.

Several works  have analyzed the complexity of this method~\cite{flaxman2005online,agarwal2010optimal,Jamieson12}, and the upper and lower bounds strongly depend on the dimension of the search space. The upper bounds also typically depend on the largest magnitude reward $B$.  If the function values are noisy, even for convex functions, the convergence rate is $O((d^2B^2/T)^{-1/3})$, and this assumes you get the algorithm parameters exactly right. For strongly convex functions, this can be reduced to $O((d^2B^2/T)^{-1/2})$ function evaluations, but this result is also rather fragile to the choice of parameters. Finally, note that just adding an constant offset to the reward dramatically slows down the algorithm. If you start with a reward function whose values are in $[0,1]$ and subtract one million from each reward, this will increase the running time of the algorithm by a factor of a million, even though the ordering of the rewards amongst parameter values remains the same.

\section{Numerical comparisons}

The preceding analyses of the RL paradigms when applied to LQR are striking. A model-based approach combining supervised learning and robust control achieves nearly optimal performance given its sampling budget. Approximate dynamic programming appears to fare worse in terms of worst-case performance. And direct policy search seems to be of too high variance to work in practice. In this section, we implement these various methods and test them on some simple LQR instances to see how these theoretical predictions reflect practice.

\subsection{A Double Integrator}

As a simple test case, consider the classic problem of a discrete-time double integrator with the dynamical model
\begin{equation}\label{eq:double-int-dynamics}
    x_{t+1} = \begin{bmatrix} 1 & 1 \\ 0 & 1 \end{bmatrix} x_t
    + \begin{bmatrix} 0 \\ 1 \end{bmatrix} u_t
\end{equation}
Such a system could model, say, the position (first state) and velocity (second state) of a unit mass object under force $u$.

As an instance of LQR, we can try to steer this system to reach point $0$ from initial condition $x0 = [-1,0]$ without expending much force:
\begin{equation}\label{eq:double-int-cost}
    Q = \begin{bmatrix} 1 & 0 \\ 0 & 0 \end{bmatrix} \qquad\qquad R = r_0
\end{equation}
for some scalar $r_0$.  Note that even in this simple instance there is an element of design: Changing the value of $r_0$ will change the character of the control law balancing expending energy versus speed or reaching the desired destination.

To compare the different approaches, I ran experiments on this instance with a small amount of noise ($e_t$ zero mean with covariance $10^{-4} I$), and training episode length $L=10$. The goal was to design a controller that works on an arbitrarily long time horizon using the fewest number of simulations of length $L$.

With one simulation (10 samples) using a white noise input with unit variance, the nominal estimate is correct to 3 digits of precision. And, not surprisingly, this returns a nearly optimal control policy. Right out of the box, this nominal control strategy works well on this simple example. Note that using a least squares estimator makes the nominal controller's life hard here because all prior information about sparsity on the state transition matrices is discarded.  In a more realistic situation, the only parameter that would need to be estimated would be the $(2,1)$ entry in $B$ which governs how much force is put out by the actuator and how much mass the system has.

\begin{figure*}[t]
	      \centering
	        \includegraphics[height=1.75in]{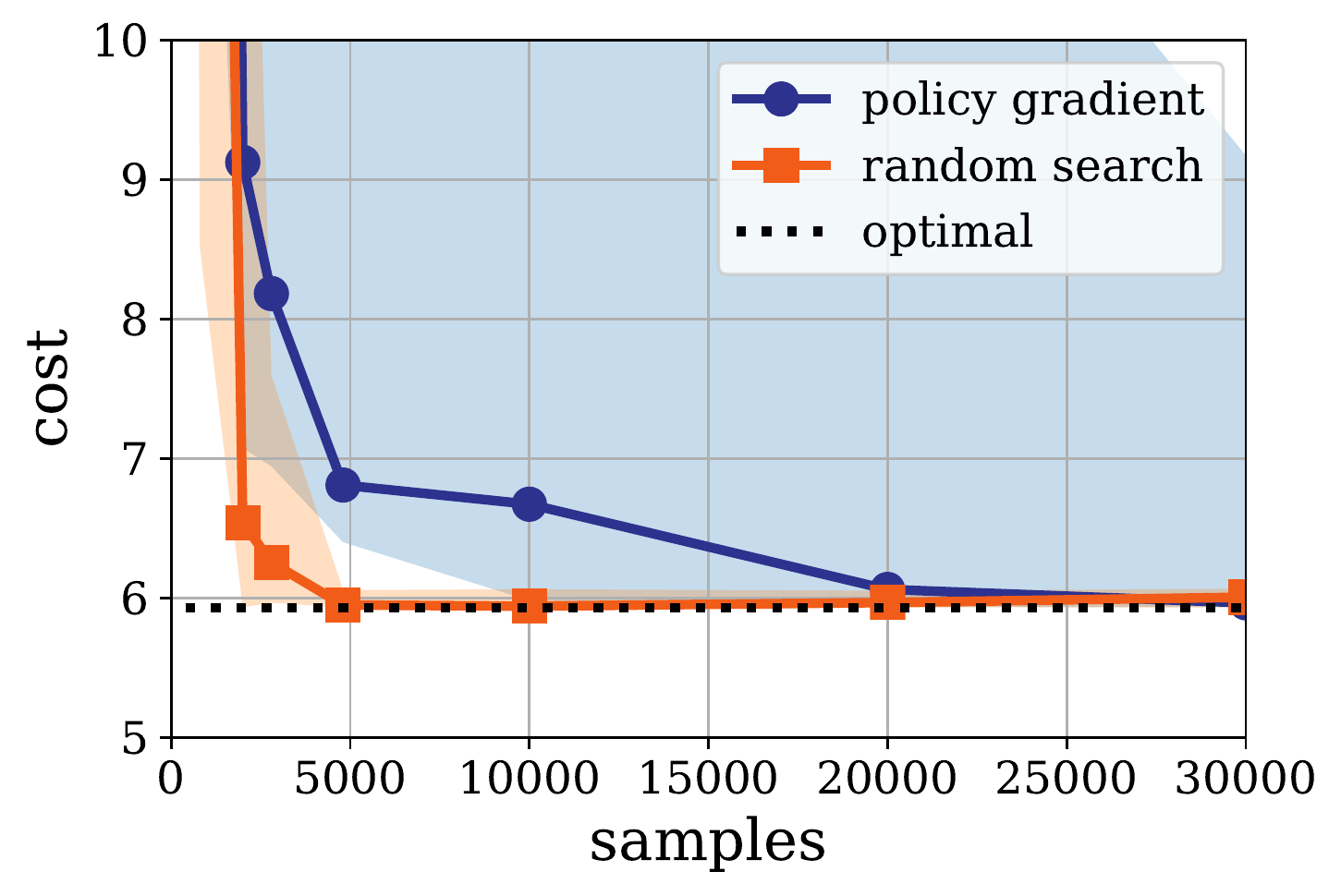}
  \caption{
Cost for the double integrator model for various RL algorithms. The solid plots denote the median performance. Shaded regions capture the maximum and minimum performance. Nominal control and LSPI are indistinguishable from the optimal controller in this experiment and hence are omitted.}\label{fig:2int-lqr}
\end{figure*}

Now, let's compare with approximate dynamic programming and policy search methods. For policy search, let us restrict to policies that use a static, linear gain as would be optimal on an infinite time horizon. Note that a static linear policy works almost as well as a time-varying policy for this simple LQR problem with two state dimensions. Moreover, there are only \emph{two decision variables} for this simple problem. For Policy Gradient, I used the Adam algorithm to shape the iterates~\cite{kingma2014adam}. I also subtracted the mean reward of previous iterates, a popular \emph{baseline subtraction heuristic} to reduce variance (Dayan~\cite{dayan1991reinforcement} attributes this heuristic to Sutton~\cite{sutton1984temporal} and Williams~\cite{Williams88}).  I was unable to get policy gradient to converge without these additional algorithmic ornamentations. I also compared against a simple ADP method based called Least Squares Policy Iteration proposed by Lagoudakis and Parr~\cite{lagoudakis2003least}. I ran each of these methods using 10 different random seeds. \fig{2int-lqr} plots the median performance of the various methods with error bars  encompassing the maximum and minimum over all trials.  Both nominal control and LSPI are able to find high quality controllers with only ten observations. Direct policy methods, on the other hand, require many times as many samples. Policy gradient, in particular requires thousands of times as many samples as simple nominal control.

\subsection{Unstable Laplacian dynamics}

As an illustrative example of the power of LQR as a baseline, let's now move to a considerably harder instance of LQR and show how it highlights issues of robustness and safety.  Consider an idealized instance of ``data center cooling,'' a popularized application of reinforcement learning~\cite{gao2014machine}.

Define the model to have three heat sources coupled to their own cooling devices.  Each component of the state $x$ is the internal temperature of one each heat source, and the sources heat up under a constant load.  They also shed heat to their neighbors. This can be approximately modeled by a linear dynamical system with state-transition matrices
\begin{equation*}
A = \begin{bmatrix} 1.01 & 0.01 & 0\\ 0.01 & 1.01 & 0.01 \\ 0 & 0.01 & 1.01 \end{bmatrix}
\qquad \qquad B = I\,.
\end{equation*}

Note that the open loop system here is \emph{unstable}: With any nonzero initial condition, the state vector will blow up because the limit of $A^k$ is infinite. Moreover, if a method estimates one of the diagonal entries of $A$ to be less than $1$, we might guess that this mode is actually stable and put less effort into cooling that source. So it is imperative to obtain a high quality estimate of the system's true behavior for near optimal control. Or rather, we must be able to ascertain whether or not our current policy is safe or the consequences can be disastrous. 

Let's try to solve the LQR problem with the settings $Q = I$ and $R= 1000 I$. This models a high relative cost for power consumption and hence may encourage small control inputs on modes that are estimated as stable. What happens for our RL methods on this instance?

\begin{figure*}[t]
\centering
\includegraphics[height=1.75in]{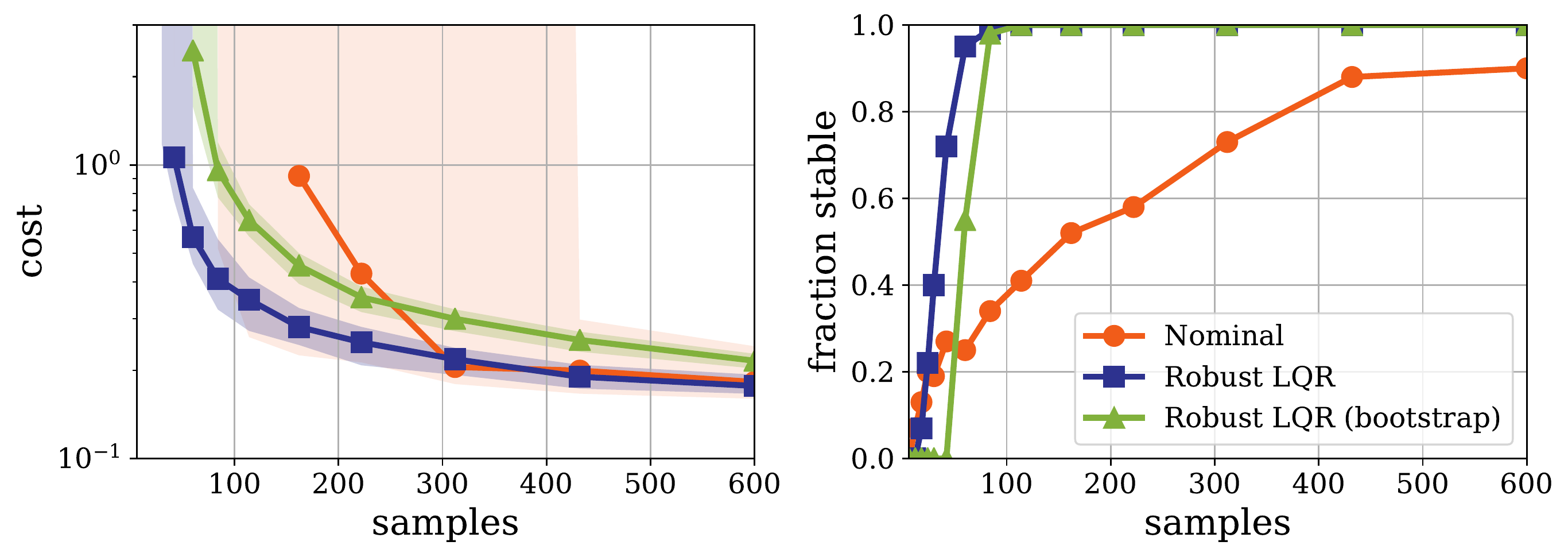} 
\caption{
(left) Cost for the Laplacian model for varied models. The blue curve shows the performance of robust LQR when provided with the true distance between the estimate and model.  The green curve displays the performance when the uncertainty is learned from data. (right) the fraction of the time that the synthesized control strategy returns a stabilizing controller.}\label{fig:datacen-lqr1}
\end{figure*}

\fig{datacen-lqr1} compares nominal control to two versions of the robust LQR problem described in section~\ref{sec:nominal}. To solve the robust LQR problem, we end up solving a small semidefinite programming problem as described by Dean et al~\cite{Dean17}. These semidefinite programs are solved on my laptop in well under a second. The blue line denotes performance when we tell the robust optimization solver what the actual distance is from the nominal model to the true model. The green curve depicts what happens when we estimate this difference between the models using a bootstrap simulation~\cite{efron79,shao2012jackknife}. Note that estimating the error from data only yields slightly worse LQR performance than exactly knowing the true model error.

Note also that the nominal controller does tend to frequently find controllers that fail to stabilize the true system. A necessary and sufficient condition for stabilization is for the matrix $A+BK$ to have all of its eigenvalues to be less than 1. We can plot how frequently the various search methods find stabilizing control policies when looking at a finite horizon in \fig{2int-lqr} (right). The robust optimization really helps here to provide controllers that are guaranteed to find a stabilizing solution. On the other hand, in industrial practice nominal control does seem to work quite well. A great open problem is to find reasonable assumptions under which the nominal controller is stabilizing. 

\begin{figure*}[t]
\centering
\includegraphics[height=2.11in]{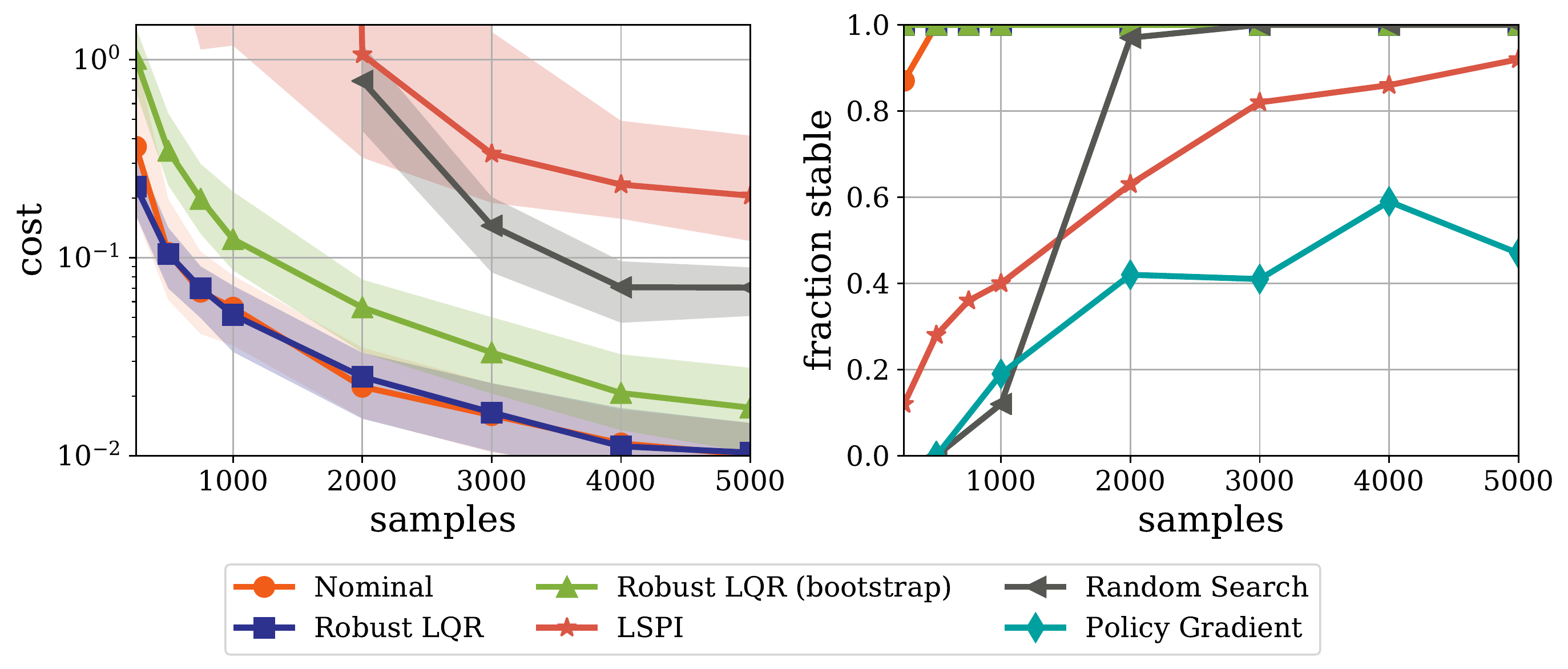}
\caption{
(left) Cost for the Laplacian model for varied models over 5000 iterations. (right) The fraction of the time that the synthesized control strategy returns a stabilizing controller.}\label{fig:datacen-lqr2}
\end{figure*}

\fig{datacen-lqr2} additionally compares the performance to model-free methods on this instance. Here we again see that they are indeed far off from their model-based counterparts. The $x$ axis has increased by a factor of 10, and yet even the approximate dynamic programming approach does not find a decent solution. Surprisingly, LSPI, which worked very well on the double integrator, now  performs worse than random search. This is likely because the LSPI subroutine requires a stabilizing controller for all iterations and also requires careful tuning of the discount factor. Not only are model-free methods sample hungry, but they fail to be safe. And safety is much more critical than sample complexity.

\section{Beyond the Linear Quadratic Regulator}

Studying simple baselines such as LQR often provides insights into how to approach more challenging problems. In this section, we explore some directions inspired by our analysis of LQR.

\subsection{Derivative Free Methods for Optimal Control}

Random search works well on simple linear problems and appears better than more complex methods like policy gradient. Does simple random search work less well on more difficult problems?

The answer, it turns out, is yes. The deep RL community has recently been using a suite of benchmarks to compare methods, maintained by OpenAI\footnote{\url{https://gym.openai.com/envs/\#mujoco}} and based on the MuJoCo simulator~\cite{todorov2012mujoco}. Here, the optimal control problem is to get the simulation of a legged robot to walk as far and quickly as possible in one direction. Some of the tasks are very simple, but some are quite difficult like the complicated humanoid models with 22 degrees of freedom. The dynamics of legged robots are well specified by Lagrange's equations~\cite{murray2017mathematical}, but planning locomotion from these models is challenging because it is not clear how to best design the objective function and because the model is piecewise linear. The model changes whenever part of the robot comes into contact with a solid object, and hence a normal force is introduced that was not previously acting upon the robot. Hence, getting robots to work without having to deal with complicated nonconvex nonlinear models seems like a solid and interesting challenge for the RL paradigm. Moreover, seminal work by Tedrake, Zhang, and Seung demonstrated that direct policy search could rapidly find feedback control policies for certain constrained legged robot designs~\cite{tedrake2004stochastic}.

Levine and Koltun were among the first to use MuJoCo as a testbed for learning-based control, and were able to achieve walking in complex simulators without special-purpose techniques~\cite{levine2013guided}. Since then, these techniques have become standard continuous control benchmarks for reinforcement learning (see, for example~\cite{silver2014deterministic, lillicrap2015continuous, schulman2015trust, schulman2015high, wu2017scalable}). Recently, Salimans and his collaborators at OpenAI showed that random search worked quite well on these benchmarks~\cite{salimans2017evolution}. In particular, they fit neural network controllers using random search with a few algorithmic enhancements. Random Search had indeed enjoyed significant success in some corners of the robotics community, and others had noted that in their applications, random search outperformed policy gradient~\cite{Stulp13}.  In another piece of great work, Rajeswaran \emph{et al} showed that Natural Policy Gradient could learn \emph{linear} policies that could complete these benchmarks~\cite{rajeswaran2017towards}. That is, they showed that static linear state feedback, like the kind we use in LQR, was also sufficient to control these complex robotic simulators.  This of course left an open question: Can simple random search find linear controllers for these MuJoCo tasks?

Guy, Mania, and I tested this out, coding up a rather simple version of random search with a couple of small algorithmic enhancements. Many RL papers were using statistics of the states and whitening the states before passing them into the neural net mapping from state to action. We found that when random search performed the same whitening with linear controls, this algorithm was able to get state-of-the-art results on all of the MuJoCo benchmark tasks~\cite{Mania18b}.

There are a few of important takeaways from this study. On the one hand, the results suggest that these MuJoCo demos are easy, or at least considerably easier than they were believed to be. Benchmarking is difficult, and having only a few simulation benchmarks encourages overfitting to these benchmarks. Indeed, it does seem like these benchmarks are more about taking advantage of simulation approximations in MuJoCo than they are about learning reasonable policies for walking. In terms of benchmarking, this is what makes LQR so attractive: LQR with unknown dynamics is a reasonable task to master as it is easy to specify new instances, and it is relatively easy to understand the limits of achievable performance.

\begin{figure*}[ht]
 \centering
 \begin{tabular}{ccc}
	      \includegraphics[height=1.5in]{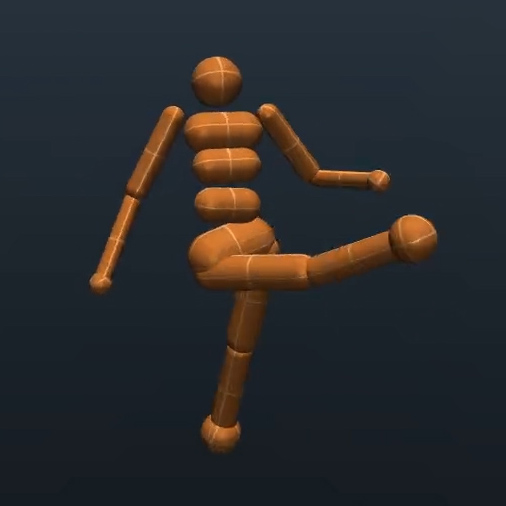} &
	      \includegraphics[height=1.5in]{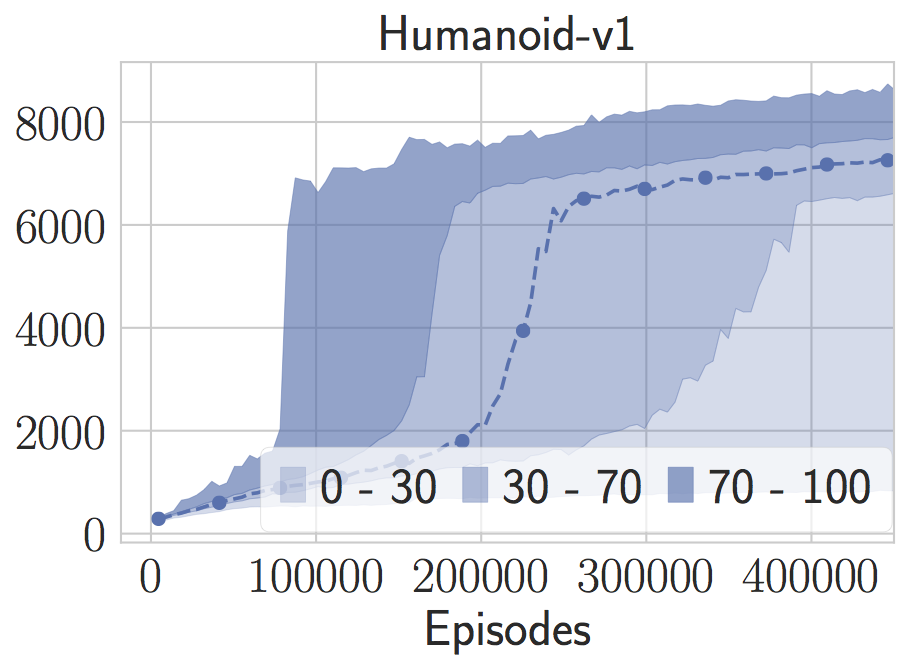} &
	      \includegraphics[height=1.5in]{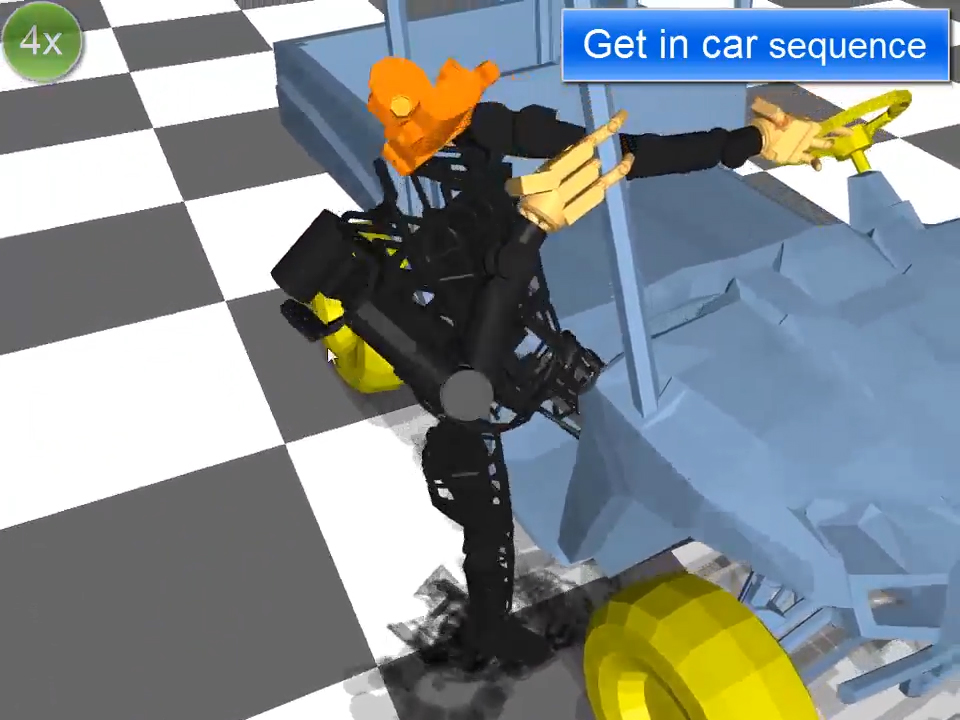} 
\end{tabular} 
  \caption{
(left) Sample frame of the MuJoCo humanoid. (middle) Variance of learning performance on 100 runs of random search on the humanoid model. Note that though high rewards are often achieved, it is more common to observe poor control performance from a random initialization. (right) Using MPC and a poor model, complex actions with humanoid simulations can be executed, such as climbing into a vehicle (image from supplementary video of Erez \emph{et al.}~\cite{erez2013integrated}).}\label{fig:mujoco}
\end{figure*}

Second, note that since our random search method is fast, we can evaluate its performance on many random seeds. All model-free methods exhibit alarmingly high variance on these benchmarks. For instance, on the humanoid task, the the model is slow to train almost a quarter of the time even when supplied with what we thought were good parameters (see Figure~\ref{fig:mujoco} (middle)). And, for those random seeds, we found the method returned rather peculiar gaits. Henderson \emph{et al} and Islam \emph{et al} observed this phenomenon with deep reinforcement learning methods, but our results on linear controllers suggest that such high variability will be a symptom of all model-free methods~\cite{henderson2017deep,islam2017reproducibility}. Though direct policy search methods are easy to code up, their reliability on any reasonable control task remains in question.

\subsection{Receding Horizon Control}

Approximate dynamic programming is closely related to canonical receding horizon control (RHC) (also known as Model Predictive Control (MPC)). In RHC an agent makes a plan based on a simulation from the present until a short time into the future. The agent then executes one step of this plan, and then, based on what it observes after taking this action, returns to short-time simulation to plan the next action. This feedback loop allows the agent to link the actual impact of its choice of action with what was simulated, and hence can correct for model mismatch, noise realizations, and other unexpected errors.

To relate RHC to ADP, note that the discounted problem
\begin{equation*}\label{eq:main-prob-disc-rhc}	
	\begin{array}{ll}	
		\mbox{maximize} &  (1-\gamma) \E_{e_t}[ \sum_{t=0}^{N-1} \gamma^t R(x_t,u_t)  + \gamma^{N}  \mathcal{Q}_\gamma(x_{N+1},u_{T+1})]\\
		\mbox{subject to} &	x_{t+1} = f(x_t, u_t, e_t)\\
		& \mbox{($x_0$ given).}
	\end{array}
\end{equation*}
is equivalent to Problem~\eq{main-prob-disc}. Here we have just unrolled the cost beyond one step. Though this is trivial, it is again incredibly powerful: the longer we make the time horizon, the less we need to worry about the Q-function being accurate. Of course, now we need to worry about the accuracy of the state-transition map, $f$. But, especially in problems with continuous variables, it is not at all obvious  which accuracy is more important in terms of finding algorithms with fast learning rates and short computation times. There is a trade-off between learning models and learning value functions, and this is a trade-off that needs to be better understood.

Though RHC methods appear fragile to model mismatch, the repeated feedback inside RHC can correct for many modeling errors. As an example, it is worth revisiting the robotic locomotion tasks inside the MuJoCo framework. These tasks were actually designed to test the power of a nonlinear RHC algorithm developed by Tassa, Erez, and Todorov~\cite{tassa2012synthesis}. The receding horizon controller works to keep the robot upright even when the model is poorly specified. Moreover, the RHC approach to humanoid control solved for the controller in 7x real time in 2012. In 2013, the same research group published a cruder version of their controller that they used during the DARPA Robotics Challenge~\cite{erez2013integrated}. All these behaviors are generated by RHC in real-time. Though the resulting walking is not of the same quality as what can be obtained from computationally intensive long-horizon trajectory optimization, it does look considerably better than the sort of gaits typically obtained by popular RL methods.
 
Is there a middle ground between expensive offline trajectory optimization and real time RHC? I think the answer is yes in the same way that there is a middle ground between learning dynamical models and learning Q-functions.  The performance of a RHC system can be improved by better modeling of the Q-function that defines the terminal cost: The better a model you make of the Q-function, the shorter a time horizon you need for simulation, and the closer you get to real-time operation.  Of course, if you had a perfect model of the  Q-function, you could just solve the Bellman equation and you would have the optimal control policy. But by having an approximation to the Q-function, high performance can still be extracted in real time.

So what if we \emph{learn} to iteratively improve the Q-function while running RHC? This idea has been explored in a project by Rosolia, Carvalho, and Borrelli~\cite{Rosolia17}. In their ``Learning MPC'' approach, the terminal cost is learned by a method akin to nearest neighbors. The terminal cost of a state is the value obtained last time that state was tried. If a state has not yet been visited, the cost is infinite. This formulation constrains the terminal condition to be in a state observed before. It enables the control system to explore new ways to decrease cost as long as it maintains the ability to reach a state that has already been demonstrated to be safe. This ``nearest-neighbors'' approach works surprisingly well in practice: in RC car demonstrations, the learned controller works better than a human operator after only a few laps around a fixed track. 

Another reason to like this blended RHC approach to learning to control is that one can hard code in constraints on controls and states and easily incorporate models of disturbance directly into the optimization problem. Some of the most challenging problems in control are how to execute safely while continuing to learn more about a system's capability, and an RHC approach provides a direct route toward balancing safety and performance.  Indeed, an interesting direction of future work would be merging the robust learning of Coarse-ID Control with receding horizon control.

\section{Challenges at the control-learning interface} 

We have set out to bridge the gap between the learning-centric views of RL and the model-centric views of control. Perhaps surprisingly, we found that for continuous control problems, machine learning seems best suited for model fitting rather than for direct control. Moreover, perhaps less surprisingly, we could seamlessly merge learned models and control action by accounting for the uncertainty in our model fits. Moreover, we showed how value functions and models could be learned in chorus and could provide impressive results on real embodied agents. These distinctions and connections are merely the beginning of what the control and machine learning communities can learn from each other. Let me close by discussing three particularly exciting and important research challenges that may be best solved with input from both perspectives.

\subsection{Merging perception and control}

One of the grand aspirations of reinforcement learning is end-to-end control, mapping directly from sensors like pixels to actions. Computer vision has made major advances by adopting an ``all-conv-net'' end-to-end approach, and many, including industrial research at NVIDIA~\cite{bojarski2016end}, suggest the same can be done for complex control tasks.

In general, this problem gets into very old intractability issues of nonlinear output feedback in control~\cite{blondel2000survey} and partially observed Markov decision processes in reinforcement learning~\cite{papadimitriou1987complexity}. Nonetheless, some early results in RL have shown promise in training optimal controllers directly from pixels~\cite{lillicrap2015continuous,mnih2015human}. Of course, these results have even worse sample complexity than the same methods trained from states, but they are making progress. 

In my opinion, the most promising approaches in this space follow the ideas of Guided Policy Search, which bootstraps standard state feedback to provide training data for a map from sensors directly to optimal action~\cite{levine2013guided,Levine16}. That is, a mapping from sensor to action can be learned iteratively by first finding the optimal action and then finding a map to that control setting. A coupling along these lines where reliance on a precise state estimator is reduced over time could potentially provide a reasonably efficient method for learning to control from sensors.

However, these problems remain daunting. Moving from fully observed scenarios to partially observed scenarios makes the  control problem exponentially more difficult. How to use diverse sensor measurements in a safe and reliable manner remains an active and increasingly important research challenge~\cite{akametalu2014reachability,aswani2013provably,berkenkamp2017safe}.

\subsection{Rethinking adaptive control}

This survey has focused on ``episodic'' reinforcement learning and has steered clear of a much harder problem: adaptive control. In the adaptive setting, we want to learn the policy online. We only get \emph{one} trajectory. The goal is, after a few steps, to have a model whose reward from here to eternity will be large. This is \emph{very} different, and \emph{much harder} that what people are doing in RL. In episodic RL, you get endless access to a simulator. In adaptive control, you get one go.

Even for LQR, the best approach to adaptive control is not settled. Pioneering work in the eighties used stochastic gradient-like techniques to find adaptive controls, but the guarantees for these schemes are all asymptotic~\cite{goodwin1981discrete}. More recently, there has been a groundswell of activity in trying to understand this problem from the perspective of \emph{online learning}. Beginning with work by Abbasi-Yadkori and Szepesvari~\cite{abbasi2011regret}, a variety of approaches have been devised to provide efficient schemes that yield near optimal control cost. Abbasi-Yadkori and Szepesvari's approach achieves an optimal reward building on techniques that give optimal algorithms for the multiarmed bandit~\cite{auer2002finite,lai1985asymptotically}. But their method requires solving a hard nonconvex optimization problem as a subroutine. Follow-up work has proposed methods using Thompson sampling~\cite{abbasi15,abeille17,17.Ouyang.LQR}, approximate dynamic programming~\cite{abbasi18}, and even coarse-ID control~\cite{Dean18}, though no method has been found that is efficiently implementable and achieves optimal cost. Again, this emphasizes that even the simple LQR problem is not at all simple. New techniques must be developed to fully understand this baseline problem, but it is clear that insights from both machine learning and control will be necessary to develop efficient adaptive control that can cope with an ever-evolving world.

\subsection{Humans in the loop}

One final important problem, which might be the most daunting of all, is how machines should learn when humans are in the loop. What can humans who are interacting with the robots do and how can we model human actions?  Though economists may have a different opinion, humans are challenging to model.

One popular approach to modeling human-robot interaction is game theoretic. Humans can be modeled as solving their own optimal control problem, and then the human's actions enter as a disturbance in the optimal control problem~\cite{sadigh2016planning}. In this way, modeling humans is similar to modeling uncertain dynamic environments. But thinking of the humans as optimizers means that their behavior is constrained. If we know the cost, then we get a complex game theoretic version of receding horizon control~\cite{bialas1989cooperative,li2017game}. But as is usually the case, humans are bad at specifying their objectives, and hence what they are optimizing must be learned. This becomes a problem of \emph{inverse optimal control}~\cite{kalman1964linear} or \emph{inverse reinforcement learning}~\cite{ng2000algorithms}, where we have to estimate the reward functions of the human and understand the loss accrued for crudely modeling these rewards.

\subsection{Towards Actionable Intelligence}

As I've expressed before, I think that all of the daunting problems in machine learning are now RL problems. Whether they be  autonomous transportation system or seemingly mundane social network engagement systems, actively interacting with reality has high stakes. Indeed, \emph{as soon as a machine learning system is unleashed in feedback with humans, that system is a reinforcement learning system.} The broad engineering community must take responsibility for the now ubiquitous machine learning systems and understand what happens when we set them loose on the world.

Solving these problems is going to need advances in both machine learning and control.  Perhaps this intersection needs a new name so that researchers can stop arguing about territory. I personally am fond of \emph{actionable intelligence} as it sums up not only robotics but smarter, safer analytics. But I don't really care what we call it: There is a large community spanning multiple disciplines that is invested in making progress on these problems. Hopefully this tour has set the stage for a lot of great research at the intersection of machine learning and control, and I'm excited to see what progress the communities can make working together.

\section*{Acknowledgements}

Countless individuals have helped to shape the contents here. First, this work was generously supported in part by two forward-looking programs at DOD: the Mathematical Data Science program at ONR and the Foundations and Limits of Learning program at DARPA. Second, this survey was distilled from a series on my blog \url{argmin.net}. I greatly appreciated the lively debates on Twitter, and I hope that even those who disagree with my perspectives here find their input incorporated into this survey. 

I'd like to thank Chris Wiggins for sharing his taxonomy on machine learning, Roy Frostig for shaping the Section~\ref{sec:policy-search}, Pavel Pravdin for consulting on how to get policy gradient methods up and running, Max Raginsky for perspectives on adaptive control and translations of Russian. I'd like to thank Moritz Hardt, Eric Jonas, and Ali Rahimi for helping to shape the language, rhetoric, and focus of the blog series and this survey. I'd also like to thank Nevena Lazic and Gergely Neu for many helpful suggestions for improving the readability and accuracy of this manuscript.  Additionally, I'd like to thank my other colleagues in machine learning and control for many helpful conversations and pointers about this material: Murat Arcak, Karl Astrom, Francesco Borrelli, John Doyle, Andy Packard, Anders Rantzer, Lorenzo Rosasco, Shankar Sastry, Yoram Singer, Csaba Szepesvari, Claire Tomlin,  and Stephen Wright.  I'd also like to thank my colleagues  in robotics, Anca Dragan, Leslie Kaebling, Sergey Levine, Pierre-Yves Oudeyer, Olivier Sigaud, Russ Tedrake, and Emo Todorov for sharing their perspectives on the sorts of RL and optimization technology works for them and the challenges they face in their research.

I'd like to thank everyone who took CS281B with me in the Spring of 2017 where I first tried to make sense of the problems in learning to control. And most importantly, a big thanks everyone in my research group who has been wrestling with these ideas with me for the past several years and for who have done much of the research that shaped my views on this space. In particular, Ross Boczar, Nick Boyd, Sarah Dean, Animesh Garg, Aurelia Guy, Qingqing Huang, Kevin Jamieson, Sanjay Krishnan, Laurent Lessard, Horia Mania, Nik Matni, Becca Roelofs, Ugo Rosolia, Ludwig Schmidt, Max Simchowitz, Stephen Tu, and Ashia Wilson.

Finally, special thanks to Camon Coffee in Berlin for letting me haunt their shop while writing.

{\small
\bibliographystyle{abbrv}
\bibliography{/Users/brecht/LaTeX/bib/brecht}
}

\end{document}